\documentclass[12pt, draft]{article} 
 \usepackage{dsfont}
 \usepackage{ams}
 \usepackage{latexsym}
   \font\twelvebm                       = cmmib10 at 12truept
   \font\tenbm                          = cmmib10 at 10truept
   \font\sevenbm                        = cmmib10 at 7truept
 = \tenbm \scriptscriptfont9              = \sevenbm
 
   at 10truept  at 10truept  at
 10truept  at 10truept
  at 10truept

 \mathchardef \BGamma            = "0900 \mathchardef \BDelta
 = "0901 \mathchardef \BTheta            = "0902 \mathchardef
 \BLambda           = "0903 \mathchardef \BXi               = "0904
 \mathchardef \BPi               = "0905 \mathchardef \BSigma
 = "0906 \mathchardef \BUpsilon          = "0907 \mathchardef \BPhi
 = "0908 \mathchardef \BPsi              = "0909 \mathchardef
 \BOmega            = "090A \mathchardef \Balpha            = "090B
 \mathchardef \Bbeta             = "090C \mathchardef \Bgamma
 = "090D \mathchardef \Bdelta            = "090E \mathchardef
 \Bepsilon          = "090F \mathchardef \Bzeta             = "0910
 \mathchardef \Beta              = "0911 \mathchardef \Btheta
 = "0912 \mathchardef \Biota             = "0913 \mathchardef
 \Bkappa            = "0914 \mathchardef \Blambda           = "0915
 \mathchardef \Bmu               = "0916 \mathchardef \Bnu
 = "0917 \mathchardef \Bxi               = "0918 \mathchardef \Bpi
 = "0919 \mathchardef \Brho              = "091A \mathchardef
 \Bsigma            = "091B \mathchardef \Btau              = "091C
 \mathchardef \Bupsilon          = "091D \mathchardef \Bphi
 = "091E \mathchardef \Bchi              = "091F \mathchardef \Bpsi
 = "0920 \mathchardef \Bomega            = "0921 \mathchardef
 \Bvarepsilon       = "0922 \mathchardef \Bvartheta         = "0923
 \mathchardef \Bvarpi            = "0924 \mathchardef \Bvarrho
 = "0925 \mathchardef \Bvarsigma         = "0926 \mathchardef
 \Bvarphi           = "0927
 \mathchardef \bA        = "0941 \mathchardef \bB        = "0942
 \mathchardef \bC        = "0943 \mathchardef \bD        = "0944
 \mathchardef \bE        = "0945 \mathchardef \bF        = "0946
 \mathchardef \bG        = "0947 \mathchardef \bH        = "0948
 \mathchardef \bI        = "0949 \mathchardef \bJ        = "094A
 \mathchardef \bK        = "094B \mathchardef \bL        = "094C
 \mathchardef \bM        = "094D \mathchardef \bN        = "094E
 \mathchardef \bO        = "094F \mathchardef \bP        = "0950
 \mathchardef \bQ        = "0951 \mathchardef \bR        = "0952
 \mathchardef \bS        = "0953 \mathchardef \bT        = "0954
 \mathchardef \bU        = "0955 \mathchardef \bV        = "0956
 \mathchardef \bW        = "0957 \mathchardef \bX        = "0958
 \mathchardef \bY        = "0959 \mathchardef \bZ        = "095A
 \mathchardef \ba        = "0961 \mathchardef \bb        = "0962
 \mathchardef \bc        = "0963 \mathchardef \bd        = "0964
 \mathchardef \bee       = "0965 
 \mathchardef \bff       = "0966 \mathchardef \bg        = "0967
 \mathchardef \bh        = "0968
 \mathchardef \bj        = "096A \mathchardef \bk        = "096B
 \mathchardef \bl        = "096C \mathchardef \bm        = "096D
 \mathchardef \bn        = "096E \mathchardef \bo        = "096F
 \mathchardef \bp        = "0970 \mathchardef \bq        = "0971
 \mathchardef \br        = "0972 \mathchardef \bs        = "0973
 \mathchardef \bt        = "0974 \mathchardef \bu        = "0975
 \mathchardef \bv        = "0976 \mathchardef \bw        = "0977
 \mathchardef \bx        = "0978 \mathchardef \by        = "0979
 \mathchardef \bz        = "097A

 \font\tencb            = cmssbx10 scaled \magstep4 \font\eigcb
 = cmssbx10 scaled \magstep2 \textfont8             = \tencb
 \mathchardef\bAs       = "1841
 \def\Asem#1#2{\mathop{\vrule height10.5pt depth5.5pt width0pt\bAs}_{#1}^{#2}}
 \def\asem#1#2{
          \ifmmode
         \ifinner
            \raise0.9pt\hbox{$\scriptstyle\bAs$}_{#1}^{#2}
         \else
            \Asem{#1}{#2}
         \fi
          \fi
          }

 \newtheorem{theo}{\small\bf Theorem}[section]
 \newtheorem{lem}{\small\bf Lemma}[section]
 \newtheorem{prop}{\small\bf Proposition}[section]
 \newtheorem{rem}{\small\bf Remark}[section]
 \newenvironment{REM}{\begin{rem} \rm}{\end{rem}}
 \newtheorem{defi}{\small\bf Definition}[section]
 \newenvironment{DEFI}{\begin{defi} \rm}{\end{defi}}
 \newtheorem{cor}{\small\bf Corollary}[section]

 \newcommand{\be}{\begin{equation}}
 \newcommand{\ee}{\end{equation}}

 \newenvironment{pr}[1]{{\small\bf {#1}:}}{}

 \newcommand{\Ca}{\widetilde{\mathds{C}}}
 \newcommand{\Ra}{\widetilde{\mathds{R}}}
 \newcommand{\Cb}{\overline{\mathds{C}}}
 \newcommand{\Rb}{\overline{\mathds{R}}}

 \newcommand{\Cinf}{\C^{\infty}}
 \newcommand{\Db}{\overline{D}}
 \newcommand{\db}{d}
 \newcommand{\dc}{\overline{d}}

 \newcommand{\chib}{\chi}
 \newcommand{\chic}{\widetilde{\chi}}

 \newcommand{\A}{A(D)}
 \newcommand{\Aa}{\widetilde{A}(D)}
 \newcommand{\Ab}{\overline{A}(D)}
 \newcommand{\Afin}{\overline{A}_{\circ}(D)}
 \newcommand{\Ainf}{\overline{A}_{\infty}(D)}

 \newcommand{\ds}{\displaystyle}
 \renewcommand{\R}{\mathds{R}}
 \newcommand{\zi}{\zeta}
 \newcommand{\Arg}{\mbox{Arg}}
 \newcommand{\Argb}{\mbox{\footnotesize Arg}}
 \newcommand{\Argc}{\mbox{\tiny Arg}}
 \renewcommand{\Im}{\mbox{Im}}

 \newcommand{\Reb}{\mbox{\footnotesize Re}}

 \newcommand{\plin}{\setminus}
 \newcommand{\h}{\eta}
 \newcommand{\dist}{\mbox{dist}}
 \newcommand{\Y}{{\cal Y}}
 \newcommand{\W}{{\cal W}}

 \title{ \Large\bf Another extension of the disc
 algebra\footnote{Work partially supported by the University of
 Athens' Research fund under Grants 70/4/4923 and
 70/4/5637.}}

 \author{\large
 V.\ Nestoridis\footnote{Corresponding author. {\tt
 e-mail address:\ vnestor@math.uoa.gr}}
 \ \ and  \ N.\ Papadatos}
 \date{\normalsize
 Department of Mathematics,
 University of Athens, \\
 Panepistemiopolis, 157 84 Athens, Greece.
 }

\begin{document}

 \maketitle

 \thispagestyle{empty}

 \begin{abstract}
 \noindent
 We identify the complex plane  $\mathds{C}$ with the open unit
 disc $D=\{z\in\mathds{C}:|z|<1\}$ by the
 homeomorphism
 $\mathds{C}\ni z \mapsto \frac{z}{1+|z|} \in
 D$. This leads to a compactification
 $\Cb$ of $\C$, homeomorphic to $\Db=\{z\in\C:|z|\leq 1\}$.
 The Euclidean  metric on $\Db$ induces a metric $\db$ on
 $\Cb$.
 We identify all uniform limits of polynomials on $\Db$
 with respect to the metric $\db$. The class of the above
 limits is an extension of the disc algebra and it is denoted
 by $\Ab$. We study properties of the elements of $\Ab$
 and topological properties of the class $\Ab$ endowed with
 its natural topology. The
 class $\Ab$ is different and,
 from the geometric point of view,
 richer
 than the class $\Aa$ introduced in
 [\ref{cite6}],  [\ref{cite7}], on the basis of the chordal
 metric $\chib$.
 \end{abstract}
 {\footnotesize {\it MSC}:  Primary 30J99, secondary 46A99, 30E10.
 \newline
 {\it Key words and phrases}: Disc Algebra; Mergelyan's Theorem;
 Polynomial Approximation.}

 \section{Introduction}
 \setcounter{equation}{0}
 \label{sec1}
 The uniform limits of the polynomials
 on the closed unit disc $\Db=\{z\in \C: |z|\leq 1\}$
 with respect to the usual Euclidean metric on
 $\C$ are exactly all functions $f:\Db\to\C$, continuous on
 $\Db$ and holomorphic in $D=\{z\in\C:|z|<1\}$.
 The class of the above functions is the disc algebra
 $\A$. In [\ref{cite6}], [\ref{cite7}], the Euclidean  metric on $\C$ is replaced by the chordal metric
 $\chib$ on the
 one-point compactification $\Ca=\C\cup\{\infty\}$ of $\C$.
 The set of uniform limits of polynomials, with respect to $\chib$,
 is an extension of $\A$ and it is denoted by $\Aa$. It
 contains the constant function equal to $\infty$ and the
 following functions: $f:\Db\to\Ca$ continuous, such that
 $f(D)\subset \C$ and $f_{|D}$ is holomorphic. The class
 $\Aa$ remains the same if we replace the chordal metric $\chib$
 by any other metric on $\Ca$ generating the same topology
 --
 the reason is that any two such metrics are uniformly
 equivalent, because $\Ca=\C\cup\{\infty\}$ is compact.
 In this sense, the set of $\chib$-uniform polynomial limits
 is an {\it invariant set}, i.e., independent of the
 specific metric one chooses for generating it.

 Instead of the one-point compactification
 $\Ca=\C\cup\{\infty\}$
 we shall consider another,
 more rich from the geometric view-point,
 compactification, as follows:
 We identify the complex plane  $\C$ with the open unit disc $D$ by the homeomorphism
 \[
 \C\ni z\mapsto \frac{z}{1+|z|}\in D.
 \]
 The natural
 compactification of $D$ is $\Db$.
 This leads to a
 compactification
 $\Cb=\C\cup\Cinf$ of $\C$,
 where the set of infinite
 points $\Cinf$ is homeomorphic to the unit circle
 $T=\partial D=\{z\in\C:|z|=1\}$; more precisely
 we write
 $\Cinf=\{\infty e^{i\theta}:\theta\in\mathds{R}\}$.
 The usual Euclidean metric on $\Db$ induces the metric
 $\db$ on $\Cb$, defined as follows:
 \be
 \label{metr}
 \hspace{-.8ex}
 \db(z_1,z_2) =
 \left\{
 \begin{array}{ll}
 \ds
 \vspace{1ex}
 \left|\frac{z_1}{1+|z_1|}-\frac{z_2}{1+|z_2|}\right|,
 &
 \mbox{for }
 z_1,z_2\in\C, \\
 \ds
 \vspace{2ex}
 \left|\frac{z_1}{1+|z_1|}-e^{i\theta}\right|,
 &
 \mbox{for }
 z_1\in\C, \ z_2=\infty e^{i\theta} \
 (\theta\in\mathds{R}),\\
 \ds
 \left|e^{i\theta_1}- e^{i\theta_2}\right|,
 &
 \mbox{for }
 z_1=\infty e^{i\theta_1}, \ z_2=\infty e^{i\theta_2} \
 (\theta_1,\theta_2\in\mathds{R}).
 \end{array}
 \right.
 \ee
 We shall investigate the uniform limits (on $\Db$)
 of the polynomials, with respect to the metric $\db$.
 The class of the limit functions is another extension
 of the disc algebra, different from $\Aa$,
 which is denoted by $\Ab$. It will be shown that it
 contains exactly two types of functions. The first type,
 the {\it finite type}, consists of all continuous
 functions $f:\Db\to\Cb$ such that $f(D)\subset\C$ and
 $f_{|D}$ is holomorphic. The second type, the {\it infinite
 type}, consists of all functions
 $f:\Db\to\Cinf\subset\Cb$ of the form
 $f(z)=\infty e^{i\theta(z)}$ where the function
 $\theta:\Db\to\mathds{R}$ is continuous on $\Db$ and
 harmonic in $D$.

 To highlight the difference between
 $\Aa$ and $\Ab$ we mention that, e.g.,
 the function $f(z)=\frac{1}{1-z}$, with $f(1)=\infty$,
 belongs to
 $\Aa$ but not to $\Ab$. An example of an element of
 $\Ab$ of the finite type, not belonging to $\A$, is
 given by $f(z)=\log\frac{1}{1-z}$ (here $f(1)=+\infty$ corresponds
 to the infinite element $\infty e^{i\theta}\in\Cinf\subset\Cb$
 with $\theta=0$).

 Furthermore, we shall investigate some properties of the
 elements of $\Ab$ and some of its topological properties
 when it is endowed with its natural metric. Finally, we
 shall consider uniform approximation with respect to the
 metric $\db$ on other compact sets, different from $\Db$.
 Certainly, due to the compactness of $\Cb$,
 all the above results remain valid if we replace $\db$ by
 any equivalent metric on $\Cb$.

 Several open questions are naturally posed and new
 directions of investigation are indicated. In particular,
 any result on $\A$ and any approximation result with
 respect to the usual Euclidean metric is worth to
 be examined in $\Ab$ and with respect to $\db$, respectively.

 \section{Preliminaries}
 \setcounter{equation}{0}
 \label{sec2}
 Let
 $\Ca=\C\cup\{\infty\}$
 be the
 one-point compactification
 of the complex
 plane $\C$,
 endowed by the chordal metric
 $\chib$, defined by
 $\chib(z_1,z_2)=\frac{|z_1-z_2|}{\sqrt{1+|z_1|^2}\sqrt{1+|z_2|^2}}$ for
 $z_1,z_2\in\C$, $\chib(z_1,\infty)=\frac{1}{\sqrt{1+|z_1|^2}}$ for
 $z_1\in\C$ and, certainly, $\chib(\infty,\infty)=0$.

 We can also define another compactification
 $\Cb$ of $\C$ with infinitely many points of infinity,
 as
 follows: First consider the homomorphism
 $G:\C\to \{w\in\C: |w|<1\}=D$ given by
 $G(z)=\frac{z}{1+|z|}$. Since $\Db$ is a compactification
 of $D$, it induces a compactification
 $\Cb$ of $\C$. The set $\Cinf$, which consists of
 all points of infinity,
 is homeomorphic with the circle $T=\partial D=\{w\in \C:
 |w|=1\}$; thus, $\Cb=\C\cup\Cinf$. Every element of $\Cinf$
 is determined by a unimodular complex number
 $e^{i\theta}$, $\theta\in\R$, and we shall denote the
 corresponding element of $\Cinf\subset \Cb$ by
 $\infty e^{i\theta}$. In particular, this compactification
 contains the usual two points of infinity of the
 real numbers ($\pm\infty$),
 in the sense that $+\infty$ corresponds to $\theta=0$ while
 $-\infty$ is related to $\theta=\pi$. The Euclidean metric
 on $\Db$ induces the metric $\db$ on $\Cb$, which is defined
 in (\ref{metr}).
 The definition of the metric $\db$ can be
 simplified if we extend $G(z)=\frac{z}{1+|z|}$, defined
 for $z\in\C$, to the points of infinity
 $z=\infty e^{i\theta}\in \Cinf$
 ($\theta\in\R$), as
 $G(z)=G(\infty e^{i\theta})=e^{i\theta}$.
 Then, $\db(z_1,z_2)=|G(z_1)-G(z_2)|$ for all $z_1,z_2\in \Cb$.

 Both metrics $\chib$ and $\db$ induce the usual topology on
 $\C$. For if $W\subset \C$ is compact, then $\chib$, $\db$ and
 the Euclidean metric
 (all restricted on $W$)
 are two by two uniformly equivalent. Thus, if $S$ is any
 set and $f,f_n:S\to W$ ($n=1,2,\ldots$) are functions,
 then the uniform convergence on $S$, $f_n\to f$,
 as $n\to+\infty$, with respect to any one of these
 metrics, implies uniform convergence with respect to the
 other two. Further, we have the following

 \begin{lem}
 \label{lem2.1}
 {\rm
 If $S$ is any set and $f,f_n:S\to \C$ ($n=1,2,\ldots$) are
 functions, then the uniform convergence
 $f_n\to f$, as $n\to\infty$, with respect to the usual
 Euclidean metric on $\C$, implies the uniform convergence
 $f_n\to f$, as $n\to\infty$, with respect to the
 metric $\db$.
 }
 \end{lem}
 \begin{pr}{Proof}
 The result is implied by the fact that
 for all $z_1,z_2\in\C$,
 \be
 \label{ineq}
 \db(z_1,z_2) \leq |z_1-z_2|.
 \ee
 This inequality can be proved as follows. Write
 \[
 \db(z_1,z_2)=\frac{|z_1-z_2+w|}{1+|z_1|+|z_2|+|z_1z_2|}, \ \
 \mbox{where} \ \
 w=z_1|z_2|-z_2|z_1|,
 \]
 and observe that
 $\ds w=\frac{|z_1|+|z_2|}{2}(z_1-z_2)-\frac{z_1+z_2}{2}(|z_1|-|z_2|)$.
 It follows that
 \[
 |z_1-z_2+w|\leq \left(1+\frac{|z_1|+|z_2|}{2} \right)
 |z_1-z_2| + \frac{|z_1|+|z_2|}{2}
 \Big||z_1|-|z_2|\Big|\leq (1+|z_1|+|z_2|)|z_1-z_2|,
 \]
 and thus,
 \[
 \db(z_1,z_2)\leq
 \left(\frac{1+|z_1|+|z_2|}{1+|z_1|+|z_2|+|z_1z_2|}\right)|z_1-z_2|\leq |z_1-z_2|.
 \]
 $\Box$
 \vspace{1em}
 \end{pr}

 Consider now the map $\Phi:(\Cb,\db)\to (\Ca,\chib)$,
 defined for all $z\in\Cb=\C\cup\Cinf$ by
 \be
 \label{map}
 \Phi(z)=\left\{
 \begin{array}{ll}
 z, & \mbox{if } z\in \C, \\
 \infty, & \mbox{if } z=\infty e^{i\theta}\in \Cinf  \ (\theta\in\R).
 \end{array}
 \right.
 \ee
 One can easily see that $\Phi$ is continuous and, therefore,
 uniformly continuous, because $\Cb$ is compact.
 This immediately implies the following
 \begin{lem}
 \label{lem2.2}
 {\rm
 If $S$ is any set and $f,f_n:S\to\Cb$ $(n=1,2,\ldots)$
 are functions, then the uniform convergence (on $S$)
 $f_n\to f$, as $n\to+\infty$, with respect to the metric
 $\db$, implies the uniform convergence
 $\Phi\circ f_n \to \Phi \circ f$, as $n\to+\infty$, with
 respect to the metric $\chib$.
 }
 \end{lem}
 \begin{rem}
 \label{rem2.1}
 {\rm
 Alternatively, one can verify Lemma \ref{lem2.2}
 by making use of the easily proved
 inequality $\chib(\Phi(z_1),\Phi(z_2))\leq 2\db(z_1,z_2)$,
 which is valid for all
 $z_1,z_2\in\Cb$.
 }
 \end{rem}
 \begin{cor}
 \label{cor2.3}
 {\rm
 If $S$ is any set and $f,f_n:S\to\C$ $(n=1,2,\ldots)$
 are functions, then the uniform convergence (on $S$)
 $f_n\to f$, as $n\to+\infty$, with respect to the usual
 Euclidean  metric on $\C$ implies the uniform convergence
 $f_n \to f$, as $n\to+\infty$, with
 respect to the metric $\chib$.
 }
 \end{cor}
 \begin{pr}{Proof} It suffices to combine Lemmas
 \ref{lem2.1} and \ref{lem2.2}, or to observe
 the trivial inequality $\chib(z_1,z_2)\leq |z_1-z_2|$,
 $z_1,z_2\in\C$.
 $\Box$
 \vspace{1em}
 \end{pr}

 Another useful fact that will be used in the sequel is the
 following; its simple proof is omitted.
 \begin{lem}
 \label{lem2.4}
 {\rm
 Let $R>0$ be a positive real number and define
 the map
 $\Phi_R:\Cb\to\C$ by
 \[
 \Phi_R(z)=\left\{
 \begin{array}{ll}
 z, & \mbox{if }
 \ds
 \vspace{1ex}
 z\in\C \ \mbox{ and } \ |z|< R, \\
 \ds
 \vspace{1ex}
 R\frac{z}{|z|}, & \mbox{if } z\in\C \ \mbox{ and } \ |z|\geq R, \\
 \ds
 \vspace{1ex}
 Re^{i\theta}, & \mbox{if } z=\infty e^{i\theta}\in\Cinf \ (\theta\in \R).
 \end{array}
 \right.
 \]
 Then, $\Phi_R$ is (uniformly) continuous.
 }
 \end{lem}

 \section{The definition}
 \setcounter{equation}{0}
 \label{sec3}
 In this section we consider the closed unit disc
 $S=\Db=\{z\in\C:|z|\leq 1\}$ and we identify the set of
 uniform limits on $\Db$ of polynomial functions with
 respect to the metric $\db$.

 Suppose $f_n$, $n=1,2,\ldots$,
 is a sequence of complex polynomials. Let $f:\Db\to \Cb$
 be a function, where $\Cb=\C\cup\Cinf$, endowed with the
 metric $\db$, is the compactification
 of $\C$ introduced previously. We assume that the sequence
 $f_n$, $n=1,2,\ldots$, converges uniformly on $\Db$
 towards $f$ with respect to the metric $\db$.
 Since polynomial functions are continuous and uniform
 convergence preserves continuity, it follows that the
 limiting function $f:\Db\to\Cb$ has to be continuous.
 Furthermore, according to Lemma \ref{lem2.2}, the sequence
 $\Phi\circ f_n=f_n$, $n=1,2,\ldots$, converges uniformly
 (with resect to $\chib$)
 to the
 function $\Phi\circ f:\Db\to\Ca=\C\cup\{\infty\}$. Here
 $\Phi$ is the map defined by (\ref{map}) and $\chib$ is the
 chordal metric on $\Ca=\C\cup\{\infty\}$
 (see section \ref{sec2}). It follows that the function
 $\Phi\circ f$ belongs to the class $\Aa$ introduced in
 [\ref{cite6}], [\ref{cite7}]. Thus, according to the definition
 of $\Aa$, the function $\Phi\circ f$ can be of the
 following two types:

 The first type contains the holomorphic functions $\Phi\circ
 f:D\to\C$,
 such that for every boundary point
 $\zi\in\partial D=T=\{w\in\C:|w|=1\}$, the limit
 \[
 \lim\limits_{{z\to\zi}, \ {z\in D}}\Phi(f(z))
 \]
 exists in
 $\Ca=\C\cup\{\infty\}$. In this case we conclude that the
 continuous function $f:\Db\to\Cb$ satisfies
 $f(D)\subset \C$ and $f_{|D}$ is holomorphic. This
 is the first type, the {\it finite type} of $\Ab$.

 The second type of elements of $\Aa$ is the constant
 function
 $\Phi\circ f\equiv\infty$; this means that for every $z\in\Db$, the
 value $f(z)$ is a point of infinity of $\Cb$.
 In this case the continuous function
 $f:\Db\to\Cb$ satisfies $f(\Db)\subset\Cinf$. Thus, $f$
 is of the form
 $f(z)=\infty e^{i\theta(z)}$ for some function
 $\theta:\Db\to\R$.
 Observe that the continuity of the function
 $f(z)=\infty e^{i\theta(z)}$, with respect to the metric
 $\db$, is equivalent to the continuity of the function
 $e^{i\theta(z)}$, with respect to the usual Euclidean  metric. Thus, the
 function
 \[
 \Db\ni z\mapsto e^{i\theta(z)}\in T=\partial D
 \]
 has to be continuous.
 Since $\Db$ is simply connected, it follows from Theorem
 5.1 of [\ref{cite5}] (see page 128)
 that the real function $\Db\ni z\mapsto \theta(z)\in \R$
 can be chosen to be continuous. Now we shall show that any
 continuous version of the function
 $\Db\ni z\mapsto \theta(z)\in \R$ is, in fact, harmonic in
 the open unit disc $D$.

 Indeed, the uniform convergence to the $\Cinf$-valued
 function $f$ (with respect to the metric $\db$) shows
 that $|f_n(z)|\to+\infty$, as $n\to+\infty$,
 uniformly on $\Db$. Thus, $|f_n(z)|\geq 1$ for all $z\in\Db$
 and for all $n\geq n_0$.
 Considering in $D$ a branch of $\log f_n$ we conclude that
 $\Arg f_n(z)=\Im[\log f_n(z)]$ is a harmonic function
 in $D$.
 We also have
 $\frac{f_n(z)}{|f_n(z)|}\to e^{i\theta(z)}$, uniformly on
 $\Db$ (with respect to the usual Euclidean metric). Thus, the same holds on $D\subset\Db$.
 It follows that
 $e^{i[\theta(z)-{\Argb} f_n(z)]}\to 1$, as $n\to+\infty$,
 uniformly on $D$. Thus, there exists $n_1\geq n_0$, such
 that for every natural number $n\geq n_1$, there exists an
 integer
 $k_n=k_n(z)\in\mathds{Z}$ such that the function
 $w_n(z):=\theta(z)-\Arg f_n(z)-2k_n(z)\pi$
 takes values in $\left(-\frac{\pi}{2},\frac{\pi}{2}\right)$
 for all $z\in D$.
 Since $e^{i w_n(z)}=\frac{e^{i\theta(z)}}{e^{i\Argc f_n(z)}}$
 is continuous in $D$ and
 $w_n(z)\in\left(-\frac{\pi}{2},\frac{\pi}{2}\right)$, it
 follows that the function $D\ni z\mapsto w_n(z)\in\R$
 is continuous ($n\geq n_1$). Therefore, writing
 $k_n(z)=\frac{1}{2\pi}(\theta(z)-w_n(z)-\Arg f_n(z))$
 we see that the function
 $D\ni z\mapsto k_n(z)\in\mathds{Z}$ is also continuous,
 and hence, constant. Thus, for $z\in D$ we may write
 $k_n(z)\equiv k_n\in\mathds{Z}$, independent of $z\in D$
 (for $n\geq n_1$). Now, the
 uniform convergence $e^{i w_n(z)}\to 1$, as $n\to+\infty$
 (with respect to the Euclidean metric) and the fact
 that $w_n(z)\in\left(-\frac{\pi}{2},\frac{\pi}{2}\right)$
 for all $z\in D$ and for all $n\geq n_1$, imply that
 $w_n(z)\to 0$, uniformly on $D$
 (with respect to the Euclidean metric).
 Equivalently, $2k_n\pi +\Arg f_n(z) \to \theta(z)$, as
 $n\to+\infty$, uniformly on $D$. Since the functions
 $2k_n\pi +\Arg f_n(z)$, $n\geq n_1$, are harmonic in $D$,
 it follows that the function $\theta(z)$, being
 a uniform limit of harmonic functions, is also
 harmonic in $D$.

 Therefore, we have shown that the second type of the limiting
 functions are continuous functions
 $f:\Db\to\Cinf\subset\Cb$, of the form
 $f(z)=\infty e^{i\theta(z)}$, where the function
 $\Db\ni z\mapsto \theta(z)\in\R$ can be chosen
 to be continuous on $\Db$ and harmonic in $D$.
 The above $\Cinf$-valued functions generate the
 {\it infinite type} of $\Ab$.

 Conversely, we shall show that each continuous function
 $f:\Db\to\Cb$ of
 the previous two types is, indeed,
 the uniform
 limit (with respect to the metric $\db$)
 of a sequence of polynomials $f_n:\Db\to\C$. The details
 are as follows:

 Suppose first that $f$ is of the finite type. That is,
 $f:\Db\to\Cb$ is continuous (with respect to the metric
 $\db$), $f(D)\subset \C$ and $f_{|D}$ is holomorphic.
 Since $\Db$ is compact, it follows that $f$ is uniformly
 continuous. Therefore, for a given $\epsilon>0$, we can
 find a real number $r=r(\epsilon)$, $0<r<1$, such that
 $\db(f(z),f(rz))<\frac{\epsilon}{2}$ for all $z\in\Db$.
 Since the function $f:D\to\C$ is holomorphic, a partial sum
 $P(w)$ of the Taylor development of $f(w)$ satisfies
 $\left|f(w)-P(w)\right|<\frac{\epsilon}{2}$ for all
 $w:|w|\leq r$. Thus,
 $\left| f(rz)-P(rz)\right|<\frac{\epsilon}{2}$
 for all $z\in\Db$. Using inequality
 (\ref{ineq}) we see that $\db(f(rz),P(rz))<\frac{\epsilon}{2}$
 for all $z\in\Db$. Now, the triangle inequality implies
 that  $\db(f(z),Q(z))<\epsilon$, $z\in\Db$, where
 $Q=Q_{\epsilon}$ is the polynomial defined by
 $Q(z)=P(rz)$. Setting $\epsilon=1/n$ we conclude
 that the sequence of polynomials $f_n=Q_{1/n}$,
 $n=1,2,\ldots$, approximates $f$ uniformly on $\Db$,
 with respect to $\db$.

 Finally suppose that $f$ is of the infinite type. That is,
 assume that $f:\Db\to\Cinf\subset\Cb$ is continuous (with
 respect to $\db$) and of the form
 $f(z)=\infty e^{i\theta(z)}$, $z\in\Db$,
 where the real-valued function $\theta(z)$, $z\in \Db$,
 is continuous on $\Db$ and harmonic in $D$. Let
 $\epsilon>0$. Since the function
 $\Db\ni z\mapsto \theta(z)\in\R$ is uniformly continuous,
 there exists $r=r(\epsilon)$, with $0<r<1$, such that
 $|\theta(z)-\theta(rz)|<\frac{\epsilon}{3}$ for all $z\in\Db$.
 Therefore, from
 $\db(f(z),f(rz))=|e^{i\theta(z)}-e^{i\theta(rz)}|
 \leq |\theta(z)-\theta(rz)|$
 we conclude that
 $\db(f(z),f(rz)) <\frac{\epsilon}{3}$
 for all $z\in\Db$.
 Since the function $D\ni z\mapsto \theta(z)\in\R$
 is harmonic in the simply connected domain $D$, we can
 find a holomorphic function $g\in {\cal H}(D)$ such that
 $\Im g(z)=\theta(z)$, $z\in D$.
 Setting $\delta=\min\{e^{\Reb g(w)}: |w|\leq r\}$
 we see that $\delta>0$.
 Thus, for any
 $n\in\{1,2,\ldots\}$ and any $w\in\C$ with $|w|\leq r$,
 we have
 \[
 \db(ne^{g(w)},f(w))=\frac{1}{1+n e^{\Reb g(w)}}\leq
 \frac{1}{1+\delta n}.
 \]
 It follows that we can fix a
 large enough $n=n(\epsilon)$ so that
 $\db(n e^{g(rz)}, f(rz))<\frac{\epsilon}{3}$ for all
 $z\in\Db$. On the other hand, since the function $n e^{g(w)}$
 is holomorphic in $D$,
 it can be approximated, uniformly over $\Db_r=\{w\in\C:|w|\leq r\}$,
 by a partial sum of its Taylor expansion, with
 respect to the Euclidean metric.
 Therefore, there exists a polynomial $P(w)$ such that
 $|n e^{g(w)}-P(w)|<\frac{\epsilon}{3}$ for all $w:|w|\leq r$.
 Using the inequality (\ref{ineq})
 (cf.\ Lemma \ref{lem2.1}) we conclude that
 $\db(n e^{g(rz)},P(rz))<\frac{\epsilon}{3}$ for all
 $z\in\Db$. By considering the polynomial
 $Q=Q_{\epsilon}$, defined by
 $Q(z)=P(rz)$, we conclude that
 for all
 $z\in \Db$,
 \[
 \db(f(z),Q_{\epsilon}(z))
 \leq
 \db(f(z),f(rz))+ \db(f(rz),n e^{g(rz)})+\db(n
 e^{g(rz)},P(rz))<
 \epsilon.
 \]
 It follows that the sequence of polynomials $f_n=Q_{1/n}$,
 $n=1,2,\ldots$, approximates $f$ uniformly on $\Db$,
 with respect to $\db$.

 Therefore, we naturally arrived  at
 the following definition.
 \begin{DEFI}
 \label{def3.1}
 Let $\Cinf=\{\infty e^{i\theta}:\theta\in\R\}$ and set
 $\Cb=\C\cup\Cinf$, endowed by the metric $\db$ defined in
 (\ref{metr}). Let also
 $D=\{z\in\C:|z|<1\}$ and $\Db=\{z\in\C:|z|\leq1\}$ be,
 as usually,
 the open and the closed unit disc, respectively. We denote
 by $\Ab$ the class of continuous functions
 $f:\Db\to\Cb$ of the following two types.

 {\small\bf (a) The finite type:} It contains the
 continuous functions $f:\Db\to\Cb$ such that
 $f(D)\subset \C$ and $f_{|D}$ is holomorphic in $D$.

 {\small\bf (b) The infinite type:} It contains the continuous
 functions $f:\Db\to\Cinf\subset\Cb$ of the form
 $f(z)=\infty e^{i\theta(z)}$, $z\in\Db$, where the
 real-valued function $\theta:\Db\to\R$ is harmonic
 in $D$ and continuous (with respect to the usual
 Euclidean metric) on $\Db$.
 \end{DEFI}

 Thus, we have shown the following

 \begin{theo}
 \label{theo3.2}
 {\rm
 Let $\Db$ be the closed unit disc. The set of
 uniform limits (on $\Db$) of the complex-valued
 polynomials -- with respect to the metric $\db$,
 defined on $\Cb=\C\cup\Cinf$ by (\ref{metr}) --
 coincides with the class
 $\Ab$.
 }
 \end{theo}

 Obviously $\Ab$ contains the disc algebra
 $\A=\{f:\Db\to \C$, continuous on $\Db$, holomorphic in
 $D\}$. It is an extension of $\A$, essentially different than
 the extension $\Aa$ obtained in [\ref{cite6}], [\ref{cite7}],
 using the chordal metric $\chib$. Indeed, it is easily seen that
 the function $f(z)=\frac{1}{1-z}$, $|z|<1$, is neither a
 restriction in $D$ of any element of $\A$, nor a
 restriction in $D$ of any element of $\Ab$. However,
 $f(z)$ can be extended to $\Db$ in an obvious manner,
 so that the resulting $\Ca$-valued function is
 $\chib$-continuous and forms an element of $\Aa$.
 Loosely speaking, $f(z)=\frac{1}{1-z}$ belongs
 to $\Aa$ but not to $\Ab$.

 Note that any conformal mapping of $D$ onto an open strip
 (or half-strip) belongs to $\Ab$ and not to $\A$.
 More generally,
 it is also true that
 for any nonzero complex numbers $c_k$, $k=1,2,\ldots,n$,
 the function
 \[
 f(z)=\sum_{k=1}^n c_k \log\frac{1}{e^{i\theta_k} -z}, \
 \  (\mbox{where } 0\leq \theta_1<\theta_2<\cdots<\theta_n<2\pi)
 \]
 is an element of $\Ab$, and not of $\A$.
 In fact, the above examples belong
 to the finite type of $\Ab$. It is easy to show that
 if $f$ belongs to $\Ab$ (and is of finite type)
 then $\Phi\circ f$ belongs to $\Aa$ -- see
 (\ref{map}), Lemma \ref{lem2.2} and
 Remark \ref{rem2.1}. On the other hand
 the converse does not hold, as we can see by the
 example $f(z)=\frac{1}{1-z}$; obviously, this function cannot
 be written as $f=\Phi\circ g$ for some $g\in\Ab$.

 Some simple examples of elements of $\Ab$ of the infinite
 type can be constructed as follows:
 Consider a polynomial $P$ not vanishing at any point of
 $\Db$. Then, the sequence $nP(z)$, $n=1,2,\ldots$,
 converges uniformly on $\Db$, with resect to the metric
 $\db$, to the function $f(z)=\infty e^{i\Argb P(z)}$
 which, certainly, belongs to the infinite type of $\Ab$.
 In particular, taking $P(z)=(2+z)^k$ with large enough
 $k\in\mathds{N}$ (in fact, $k=6$ suffices),
 we see that the image of the
 limiting function $f$ covers the whole $\Cinf$.

 \section{Some properties of the elements of $\Ab$}
 \setcounter{equation}{0}
 \label{sec4}
 Let $f\in\Ab$. Then $\Phi\circ f\in \Aa$ and, applying
 Proposition 3.1 of [\ref{cite6}], [\ref{cite7}],
 we obtain the following.
 \begin{prop}
 \label{prop4.1}
 {\rm
 Let $T=\partial D=\{\zi\in\C:|\zi|=1\}$ be the unit circle and assume that
 $f\in\Ab$.
 \\
 (a)
 If for some $c\in\C$ the set $\{\zi\in T:f(\zi)=c\}$
 has positive Lebesgue measure
 then $f$ is constant.
 \\
 (b) If the set
 $\{\zi\in T:f(\zi)\in\Cinf\}$
 has positive Lebesgue measure
 then $f$ is of the infinite type.
 }
 \end{prop}
 \begin{REM}
 \label{rem4.1}
 If $f(z)=\infty e^{i\theta(z)}$ belongs to the infinite
 type of $\Ab$ then $\theta(z)$ can be constant on a subarc
 of $T$ with strictly positive length without being
 constant on $\Db$. In fact, any continuous function
 $T\ni\zi\mapsto\theta(\zi)\in\R$ has a unique extension
 $\Db\ni z\mapsto\theta(z)\in\R$ which is continuous on
 $\Db$ and harmonic in $D$; it defines a unique $f\in\Ab$
 of infinite type.
 \end{REM}
  \begin{prop}
 \label{prop4.2}
 {\rm
 Let $K$ be a
 compact subset of $T$
 having Lebesgue measure zero. Then any
 continuous (with respect to the metric $\db$)
 function $\phi:K\to\Cinf$ is the restriction
 of some $f\in \Ab$, of finite type,
 such that $f^{-1}(\Cinf)=K$.
 }
 \end{prop}
 \begin{pr}{Proof}
 There exists a function $g\in\A$ such that
 $g(z)=1$ on $K$ and $|g(z)|<1$ on $\Db\plin K$
 (see [\ref{cite2}], p.\ 81).
 Also, since $\phi(K)\subset \Cinf$, we can write
 $\phi(\zi)=\infty e^{i\theta(\zi)}$, $\zi\in K$,
 for some real-valued function $\theta:K\to\R$. Since
 the function $\phi$ is continuous (on $K$)
 with respect to $\db$,
 the function $e^{i\theta(\zi)}$ is continuous
 (on $K$) with respect to the usual Euclidean metric.
 It follows from [\ref{cite2}], [\ref{cite1}], [\ref{cite9}] that there exists
 a function $h\in \A$ such that $h(\zi)=e^{i\theta(\zi)}$
 for $\zi\in K$.
 Now it is easy to verify that the function given by
 \[
 f(z)=\left\{
 \begin{array}{cl}
 h(z) \log \frac{1}{1-g(z)},
 &
 \mbox{if }z\in \Db\plin K,
 \\
 \phi(z), & \mbox{if }z\in K,
 \\
 \end{array}
 \right.
 \]
 has the desired properties.
 $\Box$
 \vspace{1em}
 \end{pr}

 If $E\subset T$ is compact with positive Lebesgue measure,
 then $E$ is not a compact of interpolation for $\Ab$. That
 is, there exists a continuous function $\h:E\to\Cb$ which
 does not have an extension $f\in\Ab$. Indeed,
 let $\zi_0\in E$ and let $V$ be an arc with middle point $\zi_0$
 and length less than the one half of the Lebesgue measure of $E$.
 We set $\h\equiv 0$ on $E\plin V$, $\h(\zi_0)=1$, and
 we extend $\h$ linearly on $V$.
 Assume now that for some $f\in\Ab$
 it is true that $f_{|E}\equiv \h$. Then
 $f_{|E\plin V}\equiv 0$. Since $E\plin V$ has positive
 Lebesgue measure, Proposition \ref{prop4.1}(a)
 implies that $f\equiv 0$, which contradicts
 \vspace{.7em}
 $f(\zi_0)=\h(\zi_0)=1$.

 \noindent
 {\small\bf Question 1} If $E\subset T$ is a compact set with
 Lebesgue measure zero, is it true that $E$ is a compact of
 interpolation of $\Ab$? That is, is it true that every
 continuous function $\h:E\to\Cb$ has an extension in
 $\Ab$?
 \newline
 [We refer to [\ref{cite1}], [\ref{cite9}]
 for the corresponding
 result for
 \vspace{.7em}
 $\A$.]

 One can easily see that for any compact set
 $E\subset T$, every continuous function $\h:E\to\Cinf$
 has an extension in $\Ab$ of the infinite type; this
 follows from the fact that every continuous function
 $\theta:E\to\R$ has a continuous extension on $\Db$ which
 is harmonic in $D$. Of course, this extension is unique only
 in the case $E=T$.

 Another question is as
 \vspace{.7em}
 follows:

 \noindent
 {\small\bf Question 2} Characterize the compact sets $E\subset \Db$
 having the property that every continuous function
 $\h:E\to\Cb$, with $\h(E\cap D)\subset \C$, has an
 extension in $\Ab$.
 \newline
 [We refer to [\ref{cite3}]
 for the corresponding
 result for
 \vspace{.7em}
 $\A$.]

 One can also pose questions on the nature
 of the zero set of a function $f\in\Ab$ of the finite
 type. Also, what can be said about the nature of the set
 $\{z\in\Db:f(z)=\infty e^{i0}\}=\{z\in\Db:f(z)
 =+\infty\}=f^{-1}(+\infty)$,
 when $f$ is of the infinite type?
 This is related to the zero sets (in $\Db$) of functions
 $\theta:\Db\to\R$ which are continuous on $\Db$
 and harmonic in $D$.

 The maximum principle does not hold in $\Ab$. Indeed,
 consider the polynomials $f(z)=z$ and $g(z)=2z$, which
 certainly belong to $\A\subset\Ab$. Then
 $\db(f(z),g(z))=\frac{|z|}{(1+|z|)(1+2|z|)}$.
 For $|z|=1$ we find
 $
 \db(f(z),g(z))=\frac{1}{6}
 <\frac{1}{3+2\sqrt{2}}
 =\db\left(f\left(\frac{1}{\sqrt{2}}\right),g\left(\frac{1}{\sqrt{2}}\right)\right).
 $
 However, we have the following
 \begin{prop}
 \label{prop4.3}
 {\rm
 Let $f,g\in\Ab$ and suppose that $f(\zi)=g(\zi)$
 for all $\zi\in T$. Then $f\equiv g$.
 }
 \end{prop}
 \begin{pr}{Proof}
 Consider the set
 $A=\{\zi\in T: f(\zi)\in \Cinf\}=\{\zi\in T: g(\zi)\in \Cinf\}$.
 If $A$ has positive Lebesgue measure then both $f$ and $g$
 are of infinite type. Write $f(z)=\infty e^{i\theta(z)}$
 and $g(z)=\infty e^{i\phi(z)}$ where
 $\theta,\phi:\Db\to\R$ are continuous functions on $\Db$,
 harmonic in $D$. Since $f(\zi)=g(\zi)$ for all $\zi\in T$
 we conclude that $\theta(\zi)=\phi(\zi)+2k\pi$ for all $\zi\in T$,
 where $k$ is an integer independent of $\zi\in T$. This
 implies that $\theta(z)=\phi(z)+2k\pi$ for all $z\in \Db$
 and, thus, $f\equiv g$.

 Suppose now that $A$ has Lebesgue measure zero. Then $f$
 and $g$ are both of finite type. Thus, $f(D)\subset \C$,
 $g(D)\subset \C$ and both $f$, $g$ are holomorphic in $D$.
 Therefore, the function $f-g$ is holomorphic in $D$ with zero limits on $T\plin A$. Since $T\plin A$ contains a compact
 set of positive Lebesgue measure, Privalov's Theorem
 ([\ref{cite4}], p.\ 84) implies $f\equiv g$. This completes
 the proof.
 $\Box$
 \end{pr}
 \begin{REM}
 \label{rem4.2}
 Assume that $f,g\in\Ab$ coincide on a compact set $E\subset T$
 with positive Lebesgue measure. If $f$ is of finite type
 then $g$ is also of finite type and Privalov's Theorem
 ([\ref{cite4}], p.\ 84) implies $f\equiv g$. If, however,
 $f$ and $g$ are of infinite type, it may happen
 $f\neq g$. For example, set
 $\theta(\zi)=0$ on $\{\zi\in T:\zi=e^{it}, \ 0\leq t\leq \pi\}$
 and consider two different continuous (real-valued) extensions
 $\theta_1$, $\theta_2$ on $T$. Extending $\theta_1$ and
 $\theta_2$ on $\Db$, using the Poisson kernel, we
 find that the functions $f(z)=\infty e^{i\theta_1(z)}$ and
 $g(z)=\infty e^{i\theta_2(z)}$ belong to $\Ab$,
 coincide on $\{e^{it}, 0\leq t\leq \pi\}$
 but $f\neq g$.
 \end{REM}

 One more question is the
 \vspace{.8em}
 following.

 \noindent
 {\small\bf Question 3}
 Applying Lebesgue dominated convergence Theorem, or by a direct computation, one can see that the function
 $f(z)=\log\frac{1}{1-z}$ satisfies the mean value
 property:
 \[
 f(0)=\frac{1}{2\pi} \int_{0}^{2\pi} f(e^{i\theta}) d\theta.
 \]
 Is this the case for all $f\in\Ab$, of finite
 \vspace{.8em}
 type?

 \section{Some topological properties of $\Ab$}
 \setcounter{equation}{0}
 \label{sec5}
 We recall that $\A=\{f:\Db\to\C$, continuous on $\Db$
 and holomorphic in $D\}$
 is a Banach algebra if it is endowed with the usual supremum
 norm. Furthermore,
 $\Aa=\{f:\Db\to\Ca=\C\cup\{\infty\}$,
 continuous on $\Db$, $f(D)\subset\C$, $f_{|D}$
 holomorphic$\}\cup
 \{f:f\equiv\infty\}$
 is a complete metric space
 if it is endowed with the metric
 $\chic$ given by (see [\ref{cite6}], [\ref{cite7}])
 \[
 \chic(f,g)=\sup_{|z|\leq 1} \chib(f(z),g(z)), \ \
 f,g\in \Aa.
 \]
 Finally, $\Ab$ is naturally endowed
 with the metric
 \[
 \dc(f,g)=\sup_{|z|\leq 1} \db(f(z),g(z)), \ \
 f,g\in \Ab.
 \]
 \begin{prop}
 \label{prop5.1}
 {\rm
 The metric space $(\Ab,\dc)$ is complete.
 The disc algebra $\A$ is an open and dense subset of
 $\Ab$. The relative topology of $\A$ from
 $\Ab$ coincides with the usual topology of
 $\A$.
 }
 \end{prop}
 \begin{pr}{Proof}
 Consider the set
 ${\cal B}=\{f:\Db\to\Cb\}$
 endowed with
 the metric
 \[
 \beta(f,g)=\sup_{|z|\leq 1}\db(f(z),g(z)),  \ \
 f,g\in {\cal B}.
 \]
 Since $(\Cb,\db)$ is complete,
 it follows that $({\cal B},\beta)$
 is complete. According to Theorem
 \ref{theo3.2}, $\Ab$ is the closure in
 $\cal B$ of the set of polynomials.
 Thus, $\Ab$ is a closed subset of the complete metric
 space $\cal B$. It follows that $\Ab$ is also complete.
 Let $f\in\A$; then $f(\Db)$ is a compact subset of
 $\C$ and, obviously, the compact sets
 $f(\Db)$ and $\Cinf$ are disjoint.
 Thus, $\dist(f(\Db),\Cinf)=\delta>0$.
 It is easily seen that if a function $g\in\Ab$
 satisfies $\dc(f,g)<\delta$ then
 $g(\Db)\subset \C$ and hence, $g\in\A$.
 Therefore, $\A$ is an open subset of $\Ab$.
 It is also dense, because it contains the set of
 polynomials which is dense, according to
 Theorem \ref{theo3.2}.

 Let $f,f_n\in\A$ ($n=1,2,\ldots$) and assume that
 $f_n\to f$, as $n\to+\infty$, in $\A$. From
 Lemma \ref{lem2.1} we easily see that
 $f_n\to f$, as $n\to+\infty$, in $\Ab$.
 Conversely, assume that
 $f_n\to f$, as $n\to+\infty$, in $\Ab$,
 for some functions $f,f_n\in\A$ ($n=1,2,\ldots$).
 Then, the set $f(\Db)$ is a compact subset of $\C$
 which is disjoint of the compact set $\Cinf$,
 so that $\dist(f(\Db),\Cinf)=\delta>0$.
 Thus, for $n\geq n_0$ we have
 $f_n(\Db)\subset\{w\in\Cb:
 \dist(w,f(\Db))\leq \delta/2\}=E$, say, which is a
 compact subset of $\C$. On $E$ the usual Euclidean  metric and the metric $\db$ are uniformly equivalent.
 It follows that $f_n\to f$, as $n\to+\infty$, in $\A$.
 This completes the proof.
 \vspace{.8em}
 $\Box$
 \end{pr}

 Consider now the function $\F:\Ab\to\Aa$ defined by
 $\F(f)(z)=\Phi(f(z))$, i.e.,
 \[
 \Ab\ni f \mapsto \F(f)=\Phi\circ f\in \Aa,
 \]
 (see Lemma \ref{lem2.2} and (\ref{map}) for the definition
 of the map $\Phi:\Cb\to\Ca$).
 According to Lemma \ref{lem2.2} the function $\F$ is
 continuous.
 \begin{cor}
 \label{cor5.2}
 {\rm
 The set $\F(\Ab)$ is a dense subset of $\Aa$; in fact it
 is residual.
 }
 \end{cor}
 \begin{pr}{Proof}
 Since $\F(\Ab)$ contains the set of polynomials, it
 is dense in $\Aa$ ([\ref{cite6}], [\ref{cite7}]).
 Also, $\F(\Ab)$ contains $\A$ which is open and dense
 in $\Aa$ ([\ref{cite6}], [\ref{cite7}]). Thus,
 $\F(\Ab)$ is residual in $\Aa$.
 \vspace{.7em}
 $\Box$
 \end{pr}

 Also, it is easily seen that
 the elements of finite type of $\Ab$ form
 an open dense subset, say $\Afin$, of $\Ab$.
 It follows that the elements of the infinite type
 of $\Ab$ form the closed subset
 $\Ainf=\Ab\plin\Afin$, which is
 of the first category. Since $\Ainf$ is closed,
 it is a $G_{\delta}$ set.

 In the following proposition we use some notation
 from [\ref{cite8}].

 \begin{prop}
 \label{prop5.3}
 {\rm
 Let $\h$ be any Hausdorff measure function. The set of
 all $f\in\Afin$ such that $\Lambda_{\h}(E_{f})=0$
 is dense and $G_{\delta}$ in $\Ab$, where
 $E_{f}=\{\zi\in T:f(\zi)\notin f(D)\}$.
 }
 \end{prop}
 The proof is similar to the proof of Proposition 4.3
 of [\ref{cite6}], [\ref{cite7}], the only difference being
 that one has to consider $f^{-1}(\Cinf)$
 in place of $f^{-1}(\infty)$.

 Next we define $\Y=\{f\in \Ab:f(D)\subset f(T)\}\subset\Ab$ and
 $\W=\{f\in \Ab:f(T)=\Cb\}\subset\Ab$.
 Arguments similar to those given in Proposition 4.5
 of [\ref{cite6}], [\ref{cite7}] show that $\Y$ is a
 non-empty closed subset of $\Ab$ of the first category.
 With a proof similar to the proof of Proposition
 4.6 in [\ref{cite6}], [\ref{cite7}], we can show that
 $\W$ is also a closed subset of $\Ab$ of the first
 category, but we do not know if $\W$ is non-empty.
 However, if we assume that every compact set $K\subset T$
 with zero Lebesgue measure is a compact of interpolation
 for $\Ab$, then we can show that $\W\neq \emptyset$.
 Indeed, let $K\subset T$ be a Cantor-type set with
 Lebesgue measure zero. It is well known that there exists
 a continuous surjection $\phi:K\to[0,1]$. Let
 $\Gamma:[0,1]\to\Db$ be a Peano curve with
 $\Gamma([0,1])=\Db$. Finally let
 $L:\Db\to\Cb$ be a homeomorphism. Then $L\circ\Gamma\circ\phi$
 is continuous on $K$ with
 $(L\circ\Gamma\circ\phi)(K)=\Cb$.
 Therefore, the assumption that $K$ is a compact of
 interpolation for $\Ab$ implies that there exists an $f\in\Ab$
 such that $f_{|K}=L\circ\Gamma\circ\phi$. This $f$ belongs to $\W$.

 \section{Concluding remarks and questions}
 \setcounter{equation}{0}
 \label{sec6}
 In the previous sections we considered uniform
 approximation by polynomials on the compact set
 $\Db$. However, we can also consider uniform approximation
 on other compact sets with respect to the metric $d$.
 Also, the approximating functions do not necessarily
 have to be polynomials.
 \begin{prop}
 \label{prop6.1}
 {\rm
 Let $L\subset\C$ be a compact set and let $z_0\in L^{0}$.
 We assume that for every boundary point
 $\zi\in\partial L$ the segment $[z_0,\zi]$ satisfies
 $[z_0,\zi]\plin \{\zi\}\subset L^{0}$. Then, the uniform
 limits, with respect to the metric $\db$, of polynomials
 on $L$ are exactly the functions
 $f:L\to\Cb$ of the following two types:

 (a) The first type (the finite type) contains the
 continuous functions  $f:L\to\Cb$ with
 $f(L^{0})\subset \C$ such that $f_{|L^{0}}$
 is holomorphic.

 (b) The second type (the infinite type) contains the
 continuous functions  $f:L\to\Cinf$ of the form
 $f(z)=\infty e^{i\theta(z)}$, where
 the function $\theta:L\to\R$ can be chosen to be
 continuous on $L$ and harmonic in $L^{0}$.
 }
 \end{prop}

 For the proof we may assume $z_0=0$ and imitate the proof
 for the case $L=\Db$. The difference is that when we
 approximate $f$ on the compact set $rL$, $0<r<1$, we have
 to use Runge's Theorem rather than considering the Taylor
 expansion of $f$. Since the approximation is uniform on $rL$
 with respect to the Euclidean distance on $\C$, Lemma
 \ref{lem2.1} shows that it is also uniform with respect to
 the metric $\db$.

 \begin{theo}
 \label{theo6.2}
 {\rm
 Let $f:T\to\Cb$ be any continuous function. Then, there
 exists a sequence of trigonometric polynomials converging
 to $f$ uniformly on $T$ with respect to the metric $\db$.
 }
 \end{theo}
 \begin{pr}{Proof} Let $\epsilon>0$. We wish to find a
 complex-valued trigonometric polynomial $Q=Q_{\epsilon}$ such that
 $\db(f(\zi),Q(\zi))<\epsilon$
 for all $\zi\in T$. According to Lemma
 \ref{lem2.4}, for any $R>0$ the
 composition $\Phi_R\circ f:T\to\C$ is continuous on $T$.
 Fix a large enough $R>0$ so that
 $\db(f(\zi),\Phi_R(f(\zi)))<\frac{\epsilon}{2}$ for all
 $\zi\in T$. Now, since $\Phi_R\circ f$ takes only finite
 complex values, it can be uniformly approximated on $T$ by
 a trigonometric polynomial $Q$ with respect to the usual
 Euclidean metric on $\C\cong\R^2$. According to Lemma
 \ref{lem2.1}, the approximation remains uniform
 for the metric $\db$. Thus, we have found a trigonometric
 polynomial $Q$ such that
 $\db(Q(\zi),\Phi_R(f(\zi)))<\frac{\epsilon}{2}$ for all $\zi\in T$.
 The triangle inequality now yields the desired result.
 \vspace{.8em}
 $\Box$
 \end{pr}

 We mention here that we do not know what is the set of the uniform limits of polynomials, with respect to the metric $\db$, on a circle.

 Let now $I$ be a compact segment in $\C$ or, more generally,
 a homeomorphic image of the segment $[0,1]$ in $\C$.
 Then, the uniform limits of the polynomials on $I$, with
 respect to the metric $\db$, are exactly
 all continuous functions $f:I\to\Cb$. The proof is similar
 to that of Theorem \ref{theo6.2}, with the difference that
 we make use of the classical Mergelyan's Theorem to
 approximate $\Phi_R\circ f$ by a (complex) polynomial.
 This is possible because $I^{0}=\emptyset$ and $\C\plin I$
 is connected.

 Fix now $I$ to be the compact interval $[-1,1]$.
 So far, we have seen approximations of
 $\Cb$-valued functions by complex-valued polynomials.
 However, when a continuous function $f:[-1,1]\to\Cb$
 is $\Rb$-valued, where
 $\Rb=\R\cup\{-\infty,+\infty\}\subset\Cb$, it is
 reasonable to approximate it, if possible,
 by real-valued polynomials, with respect
 to the metric $\db$. For the same reasoning,
 any $\chib$-continuous $\Ra$-valued
 function $f:[-1,1]\to\Ra\subset\Ca$,
 where $\Ra=\R\cup\{\infty\}$, should be
 approximated by real-valued polynomials with respect to
 the metric $\chib$.
 According to [\ref{cite6}], [\ref{cite7}], any
 continuous function $f:[-1,1]\to\Ca$, and hence,
 any continuous function $f:[-1,1]\to\Ra$, can be
 uniformly approximated by {\it complex-valued}
 polynomials with respect to the metric $\chib$;
 however,
 the approximating polynomials
 need not be real and, sometimes,
 they cannot be real.
 Consider, for example, the
 $\chib$-continuous
 function
 $f:[-1,1]\to \Ra$, given by
 \[
 f(x)=\left\{
 \begin{array}{cl}
 \ds
 \frac{1}{x}, & \mbox{if } x\in[-1,0)\cup(0,1],
 \vspace{1ex}
 \\
 \infty, & \mbox{if } x=0. \\
 \end{array}
 \right.
 \]
 Although there exist polynomial
 approximations for this $f$, it is easily
 seen that the approximating polynomials
 cannot be real-valued -- the above function
 is $\chib$-continuous and not $\db$-continuous.
 In fact,
 one can
 show that a function
 $f:[-1,1]\to\Ra$ can be uniformly approximated
 by {\it real-valued} polynomials, with respect to
 the metric $\chib$, if and only if it is of the form
 $f=\Phi\circ g$ for some $\db$-continuous function
 $g:[-1,1]\to\Rb$; here the map $\Phi:\Rb\to\Ra$ is the
 restriction on $\Rb\subset\Cb$
 of the map $\Phi$ defined in (\ref{map}).
 Also, it is easy to see that a function
 $f:[-1,1]\to\Rb$ can be uniformly approximated
 by {\it real-valued} polynomials, with respect to
 the metric $\db$, if and only if it is
 $\db$-continuous. In other words, uniform real
 polynomial approximations (with respect to the metric $\db$)
 can be found for a function $f:[-1,1]\to[-\infty,+\infty]$
 if and only if for each $x_0\in[-1,1]$,
 $\lim_{x\to x_0}f(x)=f(x_0)\in[-\infty,+\infty]$.
 If this is true then the same real polynomials
 approximate $\Phi\circ f$ in the $\chib$-metric.
 For example, the function $f:[-1,0)\cup(0,1]\to \R$
 with $f(x)=\frac{1}{x^2}$ (for $x\in[-1,1]$, $x\neq 0$)
 can be extended, in an obvious manner, to a
 $\db$-continuous function on $[-1,1]$ (setting
 $f(0)=+\infty$)  and to a
 $\chib$-continuous function on $[-1,1]$ (setting
 $f(0)=\infty$).
 It follows that this $f$ can be uniformly
 approximated by real polynomials with respect
 to the metric $\db$ (and, hence, also with respect to $\chib$).

 Finally, it is natural to ask about the uniform limits
 of polynomials on $L$, with respect
 to the metric $\db$, when $L$ is a
 compact subset of $\C$ with
 connected complement.
 Specifically, we have the
 \vspace{.7em}
 following


 \noindent
 {\small\bf Question 4}
 Let $L\subset \C$ be a
 compact set
 with connected complement.
 Let $f:L\to\Cb$ be a continuous function, such that
 for every component $V$ of $L^{0}$, the following holds:
 either $f(V)\subset \C$ and $f_{|V}$ is holomorphic,
 or $f(V)\subset \Cinf$ and $f$ is of the form
 $f(z)=\infty e^{i\theta(z)}$ for all $z\in V$, where
 the real-valued function $\theta$ is harmonic in $V$.
 Does there exist a sequence of (complex-valued) polynomials converging to $f$ uniformly on $L$
 with respect to the
 metric
 \vspace{.7em}
 $\db$?

 The existence of such a sequence of polynomials
 would lead to an extension of the classical
 Mergelyan's Theorem in the case of the metric
 $\db$; we refer to [\ref{cite10}]
 for the classical Mergelyan's
 Theorem.

 We notice that the converse is true and the proof is the same as the one given here for the particular case
 $L=\Db$. Indeed, in the proof of Theorem \ref{theo3.2}
 we have only used the fact that $D$ is a simply connected
 domain; this is the case for every component $V$ of
 \vspace{0.7em}
 $L^{0}$.

 {\bf Acknowledgement.} We thank D.\ Cheliotis
 for some helpful comments regarding the metric $\db$.

 {
 
 }
 \end{document}